\documentclass{amsart}

\usepackage{moreverb}
\usepackage{amsmath}
\usepackage{algorithmic}
 \usepackage[foot]{amsaddr}

\newcommand\BibTeX{{\rmfamily B\kern-.05em \textsc{i\kern-.025em b}\kern-.08em
T\kern-.1667em\lower.7ex\hbox{E}\kern-.125emX}}

\usepackage[T1]{fontenc}
\usepackage[norelsize]{algorithm2e}
\usepackage{epsfig}

        {\begin{enumerate}[#1]}{\end{enumerate}%
         \vspace{-.6\baselineskip}}

        {\begin{compactenum}[#1]}{\end{compactenum}}
        
  \newtheorem{Def}{Definition}

\begin{document}

\title{Towards a multigrid method for the minimum-cost flow problem}

\author{Alessio Quaglino}
\email{alessio.quaglino@usi.ch}
\author{Rolf Krause}
\email{rolf.krause@usi.ch}
\address{Institute of Computational Science, Universit\' a della Svizzera Italiana, Via G.~Buffi 13, 6900 Lugano, Switzerland}

\thanks{The support of the State Secretariat for Education, Research and Innovation SERI Swiss Space Office under grant 236-01 D3 is gratefully acknowledged.}

\subjclass[2010]{65M55, 90C51, 68W10.}

\date{\today}

\keywords{Multigrid methods, interior-point methods, parallel algorithms.}

\begin{abstract}
We present a first step towards a multigrid method for solving the min-cost flow problem. Specifically, we present a strategy that takes advantage of existing black-box fast iterative linear solvers, i.e. algebraic multigrid methods. We show with standard benchmarks that, while less competitive than combinatorial techniques on small problems that run on a single core, our approach scales well with problem size, complexity, and number of processors, allowing for tackling large-scale problems on modern parallel architectures. Our approach is based on combining interior-point with multigrid methods for solving the nonlinear KKT equations via Newton's method. However, the Jacobian matrix arising in the Newton iteration is indefinite and its condition number cannot be expected to be bounded. In fact, the eigenvalues of the Jacobian can both vanish and blow up near the solution, leading to a significant slow-down of the convergence speed of iterative solvers - or to the loss of convergence at all. In order to allow for the application of multigrid methods, which have been originally designed for elliptic problems, we furthermore show that the occurring Jacobian can be interpreted as the stiffness matrix of a mixed formulation of the weighted graph Laplacian of the network, whose metric depends on the slack variables and the multipliers of the inequality constraints. Together with our regularization, this allows for the application of a black-box algebraic multigrid method on the Schur-complement of the system. 
\end{abstract}

\maketitle


\section{Introduction}
This work is concerned with the solution of the minimum-cost flow (MCF) problem
\begin{equation} \label{LP}
\begin{split}
&\min_{x \in R^m} c^T x \\  
s.t.&
\begin{cases}
Ax=b, \\
g(x) \geq 0, \\
\end{cases}
\end{split}
\end{equation}
where $c \in R^m$ and $b \in R^n$ are known as the cost and flow vectors, $A$ is a $n \times m$ incidence matrix of a network with $n$ nodes and $m$ arcs, and the constraints $g(x)$ have the form:
\begin{equation} \label{g_of_x}
g(x) = \left(\begin{array}{c} x-x_l \\x_u - x\end{array}\right),
\end{equation}
where $x_l$ and $x_u$ are vectors of, respectively, lower and upper
bounds on the sought flow $x$. There is a long history concerning the
development of combinatorial techniques for the MCF problem. Several
classes of methods have appeared over the years, starting from the
seminal work by Ford and Fulkerson \cite{For58}. Among popular
approaches, we recall cycle-canceling algorithms \cite{Kle67,Gol88},
successive shortest paths \cite{Edm72,Tom71}, capacity-scaling
\cite{Edm72,Orl88}, cost-scaling \cite{Gol87,Gol93}, and the network
simplex \cite{Ahu92,Joh93}. Unfortunately, none of these methods has
been shown to be of optimal complexity \cite{Ahu92,Kov15}, although
they usually behave much better in practice than their theoretical
(worst case) estimates may predict.

Although already with the ellipsoid method a solution method with a
worst case polynimial complexity existed, the ``projective algorithm''
introduced by Karmarkar \cite{Kar84} not only provided a better
estimate (${\mathcal O}(n^{3.5})$ instead of (${\mathcal O}(n^6)$), but also
worked well in practice. The original algorithm by Karmarkar has been
reformulated later as interior point method (IPM). More elaborated
variants of interior point methods, i.e. primal-dual methods, nowadays
are also applied routinely to MCF problems
\cite{Cha07,Dai08,Del15,Dra12,Fra07,Oli05,Por04,Por08}. In the
framework of IPMs, the main computational burden consists in the
solution of a linear system of equations, involving the symmetric but
possibly indefinite Jacobian arising within Newton's method applied to
the KKT equations from the IPM.

As it turns out, there is no clear ``best method'', as the choice of the
``best'' method will depend on the considered application and the
targeted problem size. The practical performance of solution methods
for the MCF problem is usually evaluated using benchmarks, as
described in \cite{Gab89,Kov15}.

While combinatorial techniques are
very efficient on problems of moderate size that can run on a single
processor, they are in general not trivial to parallelize. Thus, in
terms of large scale problems, at a first glance IPMs seem to be more
attractive, since parallel linear solvers for large scale problems are
available. However in the case
of the IPM, the arising Jacobian matrices can become extremely
ill-conditioned in the neighborhood of the solution.
For this reason, the most commonly studied approach
is to apply a direct solver, such as the Cholesky factorization
\cite{Noc99}. Even with this choice, however, problem (\ref{LP}) needs
to be regularized before a factorization can be computed. A more
detailed discussion can be found in, e.g., in \cite{Fri12,Gon12},
where two related diagonal perturbations of the system are proposed.

While direct solvers in general might be considered to be
more robust,
their non-optimal complexity leads to similar shortcomings to those that plague the
combinatorial approaches, in particular the sub-optimal scaling with
respect to the number of arcs $m$. This fact motivated several authors
to consider Krylov subspace methods
\cite{Cha07,Fra07,Oli05,Por04,Por08} to solve the MCF
problem. However, the efficiency of this approach relies on the
spectral properties of the Jacobian matrix. Therefore, most of the works in this class focus on preconditioners for this
matrix. While using the diagonal of the Jacobian proved to be very
effective in the first iterations, the situation is more complex close
to the solution, where a spanning-tree preconditioner needs to be
applied. In \cite{Por04,Por08} the use of the conjugated gradient
method is proposed, employing a heuristic to select between the two
options. However, it remains an open question if the complexity of
this approach is optimal with respect to the problem size.

Outside the scope of the MCF problem, linear solvers of optimal
complexity, such as multigrid (MG) methods \cite{Bra07,Brigg00,Hack85}, are
traditionally applied in the context of linear elliptic operators. In extension, multigrid
methods can also be applied to constrained minimization problems, i.e. the 
minimization of convex quadratic functionals subject to box
constraints, as it occurs in contact problems \cite{Kra08}, or to
the IPM, see \cite{Dra12}, where the Schur's complement method is used in
presence of additional equality constraints \cite{Maa00}.
\
However, the  application of multigrid methods to the Jacobian matrices arising in the IPM is not straightforward.
This is caused by the particular structure of the MCF, i.e. the
combination of linear and inequality constraints with a linear
objective function, leads to  several difficulties. In this paper, we discuss these difficulties
and present easily implementable and robust solutions. We deal in particular with the following the following aspects:
\begin{enumerate}
\item As discussed in Section 3, the Jacobian of the KKT
          equations for the MCF problem is indefinite rather than
          positive definite. Fortunately, it turns out that this
          matrix is formally equivalent to the discretized operators
          of the mixed formulation of the Laplace operator in finite
          elements. As the latter is elliptic, this creates the bridge
          to the application of multigrid methods.
\item The Schur complement of the system, encoding the
          actual graph Laplacian operator of the graph, is applied to
          the Lagrangian multipliers of the equality constraints, rather
          than to the flow vector $x$; thus, box constraints on the
          solution become general inequality constraints for this
          variable, hindering the application of standard techniques
          for constrained elliptic problems, such as, e.g.,  projected Gauss-Seidel
          \cite{Kra08}.  Here, we develop and present alternatives.
\item The Schur complement is plagued by the same
          ill-conditioning as the full Jacobian, therefore even MG solvers
          developed for the graph Laplacian \cite{Del15,Liv12} fail to
          converge in the proximity of the true solution. Hence
          appropriate regularization techniques need to be
          studied. Our strategy is presented in Section 4.
\end{enumerate}

Our goal is to fill the gap between the study of linear solvers for
the IPM applied to the MCF and the optimal-complexity methods for
elliptic problems. More specifically, we decide to feed the Schur's
complement of the Jacobian of the KKT equations directly into a
black-box algebraic MG (AMG) preconditioner, such as, e.g., the
BoomerAMG algorithm from \cite{Fal02}. This in principle could extend
the class of problems that can be tackled efficiently using solution
methods developed for elliptic problems and improve the IPM
complexity for the MCF problem with respect to the number of
arcs. However, as explained above, a direct application of (A)MG will
lead to several difficulties. In order to overcome those, here we develop
appropriate regularizations of the problem that yield a lower and
upper bounds on the eigenvalues of the graph Laplacian. Our strategy
is to consider the MCF problem as the task of finding the static
equilibrium of a dynamic contact problem, by adding damping to the
system and combining techniques from the IPM and active-set methods
(ASM). While we focus on the MCF problem, we argue that these
techniques might help dealing with the ill-conditioning and
non-convexity of Jacobians in more general applications of the IPM.

We here follow the path of regularizing the problem, combining a
diagonal mass matrix arising from a pseudo-kinetic energy term added
to (\ref{LP}) and the active-set method (ASM) to handle inequality
constraints, as detailed in Section 4.

\subsection{Outline}The paper is organized as follows:
\begin{itemize}
\item In Section 2, we present the mathematical formulation of the IPM applied to the MCF problem.
\item In Section 3, we outline the basics of linear solvers applied to the MCF problem, highlighting why out-of-the-box MG methods cannot be used directly.
\item In Section 4, we discuss our modifications to the standard IPM to solve the KKT equations, which produce a positive definite Schur's complement.
\item In Section 5, we show some numerical results on standard benchmarks.
\end{itemize}

\section{Solving the MCF problem via IPM}
The MCF problem (\ref{LP}) yields the following first order optimality conditions
\begin{eqnarray} \label{KKT}
f_x := &~c - A^T y - \nabla g^T s~  &= 0, \\
f_y := &b - A x &= 0, \\
f_s := & -g(x)_i s_i &= 0, \quad \forall i=1\dots 2m, \\
&~ (g(x),s) ~& \geq 0, \label{complementarity}
\end{eqnarray}
where $y$ and $s$ are the Lagrange multipliers associated to the constraints $Ax=b$ and $g(x) \geq 0$, respectively, and where we have denoted the gradient of $g(x)$ as
\begin{equation} \label{nabla_g}
\nabla g := \nabla g(x) = \left(\begin{array}{c}  I \\ -I \end{array}\right).
\end{equation}
For solving the non-linear KKT conditions (\ref{KKT}), usually Newton's Method is used. The main computational effort in each Newton iteration is the solution of the linear system
\begin{eqnarray} \label{NewtonSys}
\left(\begin{array}{ccc}
0 & A^T & \nabla g^T \\
A & 0 & 0 \\
\nabla g &0 & S^{-1}G
\end{array}\right)
\left(\begin{array}{c}
\Delta x \\
\Delta y \\
\Delta s
\end{array}\right) &=&
\left(\begin{array}{c}
f_x  \\
f_y \\
S^{-1}f_s 
\end{array}\right), 
\end{eqnarray}
where $G$ and $S$ are a diagonal matrices having, respectively, $g(x)_i$ and $s_i$ on their main diagonals. The right-hand-side is evaluated based on the current Newton iterate, thereby measuring the actual residual of the solution. 

In order to find a reasonable step-size for the Newton correction, line search has to be performed to find $\alpha_x$ and $\alpha_s$ such that
\begin{eqnarray} \label{lineSearch}
x_{k+1} =& x_k + \alpha_x \Delta x &\geq 0, \\
s_{k+1} =& s_k + \alpha_s \Delta s &\geq 0.
\end{eqnarray}
In general, this straight forward application of Newton's method will lead to line search increments $\alpha_x$ and $\alpha_s$, which will be too small for making this strategy practical. To solve this issue, often reformulations based on interior point methods (IPM) are employed, see, e.g. \cite{Noc99}. In IPM, problem (\ref{LP}) is modified with a barrier potential
\begin{equation} \label{barrierLP}
\begin{split}
\min_{x \in R^m} c^T &x + \mu \sum_{i=0}^{2m} \log(g(x)_i) \\  
s.t.&
\begin{cases}
Ax=b, \\
g(x) \geq 0, \\
\end{cases}
\end{split}
\end{equation}
where $\mu>0$ is a centering parameter to be specified. By defining $s:=g(x)/\mu$, this formulation yields a modified third KKT optimality condition via the substitution
\begin{eqnarray} \label{modified_fs}
f_s \rightarrow & \mu + f_s.
\end{eqnarray}
Clearly, as $\mu \rightarrow 0$, these equations tend to the KKT conditions for the original problem (\ref{LP}). Among the several known strategies to combine accuracy (small $\mu$) with a small number of Newton iterations (large $\mu$), we decided to employ Mehrotra's predictor-corrector algorithm outlined in \cite{Noc99}. While this choice is obviously well suited for direct solvers, as it requires to solve the system twice with the same matrix, we also obtain a benefit with AMG methods, since the setup phase of the AMG has to be executed only once per Newton iteration.In our experiments, this proved to be a small benefit for graphs with a regular structure, e.g. graphs which would arise from discretizations of the Laplacian on structured meshes, but yielded a significant speedup on complicated graphs, such as the NETGEN data set discussed in the Section on numerical experiments.

\section{Mixed formulation and its limitations}
In order to be able apply multigrid methods, the Newton iteration (\ref{NewtonSys}) must be transformed into a positive definite system. First, let us observe that the last row can be solved for, reducing the system to
\begin{eqnarray} \label{normalEq}
\left(\begin{array}{ccc}
D & A^T \\
A & 0 
\end{array}\right)
\left(\begin{array}{c}
\Delta x \\
\Delta y 
\end{array}\right) &=&
\left(\begin{array}{c}
f_x  - \nabla g^T G^{-1}f_s\\
f_y
\end{array}\right),
\end{eqnarray}
where $D = -\nabla g^T G^{-1}S \nabla g$. This system is indefinite but exhibits a special structure that, in the context of FE, is known as the mixed formulation of the graph Laplacian induced by the metric $D$. We here show more precisely this equivalence, highlighting that the distortion of $D$ close to the solution renders the IPM challenging for linear solvers. We start with the mixed formulation of the Laplace equation
\begin{alignat*}{20} 
	\langle x - \nabla y,  \phi \rangle  &= 0   &\quad& \forall \phi \in E(\Omega), \\
	-\langle x,  \nabla \psi \rangle  &= \langle f, \psi \rangle &\quad& \forall \psi \in V(\Omega),
\end{alignat*}
where $E(\Omega)$ and $V(\Omega)$ are Hilbert spaces. In FE, subspaces $V_h \subset V$ and $E_h \subset E$ are chosen to approximate the solution. In particular, we select piecewise linear functions with, respectively, vertices and edges of the graph as degrees-of-freedom, and we look for a solution $(x,y) \in E_h \times V_h$. This results in a system of linear equations
\begin{alignat*}{20}
M_e &~x~& - K^T y &=  0, \\
-K &~x~& &=  M_v f.
\end{alignat*}
The matrix $K$ represents a gradient operator, using the fact that the gradient of a linear basis function is constant over the edges, we can write it as $K = A P$, thereby separating a purely combinatoric matrix $A$ from the weights matrix $P$, containing the contributions depending on the metric (i.e. the embedding of the graph into an ambient space). This yields 
\begin{alignat*}{20}
M_e &~x~&-P^T A^T y &= 0, \\
- A P &~x~& &= M_v f.
\end{alignat*}
With the substitution $\tilde{x} = Px$, we obtain
\begin{alignat*}{20}
P^{-T} M_e P^{-1} &~\tilde{x}~& + A^T y &= 0, \\
A &~\tilde{x}~& &= -M_v f, 
\end{alignat*}
whose matrix is identical to our system if $D := P^{-T} M_e P^{-1}$. Therefore, the Schur's complement of (\ref{normalEq}) is the actual Laplace operator of the graph embedded with the metric $D$
\begin{equation}
L=A D^{-1} A^T. 
\end{equation}
Since this operator is elliptic and the resulting matrix has a bounded bandwidth, it is well-known that the MG methods have optimal complexity. Unfortunately, from the definition of $D$ it becomes evident that there are three issues that hinder the straightforward application of linear solvers:
\begin{enumerate}
\item If any of the components of $s$ tends to zero, then at least one eigenvalue of $L$ blows up.
\item If any of the components of $g(x)$ tends to zero, then at least one eigenvalue of $L$ tends to zero.
\item If the gradient of the active constraints (including the equalities) is not full-rank, then the Schur's complement becomes singular. This would be equivalent to violating the inf-sup condition in FE, which results in a non-unique solution for the dual variable $y$.
\end{enumerate}
Note that while issue 3 is problem-specific and may not happen in general, 1 and 2 need to occur near the solution. In the following, we will propose a solution for each one of the above issues, allowing for a direct application of a black-box AMG solver to the Schur's complement $L$.

\section{Solving the three limitations}
As we will see, the main intuition behind the proposed regularization is to treat the MCF as the static equilibrium of a dynamic contact problem. While MG methods are well studied in this context \cite{Kra08}, the presence of the linear constraints $Ax=b$ introduces a coupling between bound constraints on the variable $y$, making techniques such as projected Gauss-Seidel \cite{Kra08} inapplicable. To see this, recall the fact that projected GS needs to project back to the feasible region one component of the solution at the time. However, box constraints on $x$ induce general inequality constraints on $y$
\begin{equation}
g(D^{-1}A^T y) \geq 0.
\end{equation}
Clearly, given a $y$ that violates this constraint, it is not trivial nor efficient to correct it while satisfying simultaneously the linear equalities. Therefore, we will introduce appropriate boundary conditions on the variable $x$ applied to the system before the Schur's complement is computed (Section 4.2). These, in turn, will make the system singular and therefore will require additional boundary conditions on the dual variable $y$ in order to remove the kernel from the system (Section 4.3). We argue that, while the dual boundary conditions are very specific to the MCF problem, the first two ideas can be applied as well to general nonlinear problems exhibiting a saddle-point structure with box constraints.

\subsection{Solving 1: Pseudo-dynamic minimization}
The first issue stems from the fact that the Hessian of (\ref{barrierLP}) with respect to $x$ vanishes. Other authors have recognized this problem and have proposed a similar regularization \cite{Fri12,Gon12}. However, we here give a physical interpretation of the perturbed problem as the problem of finding the static equilibrium of a dynamic problem, opening up the possibility to tune the method for a faster convergence through an appropriate amount of damping added to the system. We start by transforming the functional in (\ref{barrierLP}) to the Lagrangian
\begin{equation}
L(x,y) := \frac{1}{2} \dot{x}^T M \dot{x} - \mu \sum_{i=1}^{2m} \log(g(x)_i) - c^T x - y^T (Ax-b),
\end{equation}
where $M$ is an arbitrary mass matrix. In practice, we found that the choice $M = \rho I$ works reasonably well. By using a damping parameter $\beta$ and the explicit linear approximation in time $v^k = (x^k - x^{k-1}) / \Delta t$, we define an unconstrained step as
\begin{equation} \label{uncStep}
x_0^{k+1} =  x^k + \Delta t ~ {v^k} - \Delta t^2 (M^{-1} c + \beta v_k),
\end{equation}
which will serve as an initial guess of the solution $x^{k+1}$ after the Newton step. By doing so, we can use a step-and-project approach \cite{Gol07} to project back the unconstrained solution to the manifold defined by our constraints. The projection problem consists of finding the stationary point $(x,y) \in R^n \times R^m$ of the functional
\begin{equation} \label{QP}
\frac{1}{2}  \left(x - x_0^{k+1} \right) ^T M \left(x - x_0^{k+1} \right) - \Delta t^2 \left( \mu \sum_{i=1}^{2m} \log(x_i) + y^T ( Ax-b) \right).
\end{equation}
By computing the stationary equations, we obtain a modified first KKT condition
\begin{eqnarray}
\tilde{f}_x := -\frac{M}{\Delta t^2} \left( x -x_0^{k+1} \right) - A^T y - \nabla g^T s  = 0,
\end{eqnarray}
which results in the regularized version of (\ref{normalEq}):
\begin{eqnarray} \label{regularizedNormalEq}
\left(\begin{array}{cc}
D+\frac{M}{\Delta t^2} & A^T \\
A & 0  
\end{array}\right)
\left(\begin{array}{c}
\Delta x \\
\Delta y 
\end{array}\right)  &=&
\left(\begin{array}{c}
\tilde{f}_x  - \nabla g^T G^{-1} f_s  \\
f_y 
\end{array}\right).
\end{eqnarray}
It is clear to see that the Schur's complement is
\begin{equation}
\tilde{L} = A \left(\frac{M}{\Delta t^2} + D \right)^{-1} A^T,
\end{equation}
whose eigenvalues are bounded below by $M \Delta t^2$ as $s$ approaches zero. It can also be noted that for the special choice $\beta=\Delta t=1$, it follows that $\tilde{f}_x = f_x - \frac{M}{\Delta t^2} \Delta x$, so the respective right-hand sides of (\ref{normalEq}) and (\ref{regularizedNormalEq}) would coincide. This would be equivalent to the common choice of only modifying the matrix \cite{Fri12}. However, it is clear that there is a trade-off at play: bigger $\rho$ and smaller $\Delta t$ imply a better spectrum for $\tilde{L}$, but also more time steps to reach convergence. 

\begin{algorithm}
\caption{Adaptive time step computation.}
\label{adaptiveDt}
\begin{algorithmic}
\IF{$\min(s) \leq \eta$}
\STATE{$\Delta t = \sqrt{\frac{\rho}{\eta-\min(s)}}$}
\STATE{$\beta = 1$}
\ELSE
\STATE{$\Delta t = 1$}
\STATE{$\beta =0$}
\ENDIF
\end{algorithmic}
\end{algorithm}

In order to improve convergence, we performed the adaptive adjustment shown in Algorithm \ref{adaptiveDt} with the choice $\eta = 10^{-4}$, aimed at providing a lower bound on the eigenvalues of $L$ close to the solution, while taking larger steps during the first interations. We did not study further the problem of minimizing the total number of linear solver iterations required by all Newton steps. In fact, we found that for all benchmarks except the GOTO data set, we do not need to use Algorithm \ref{adaptiveDt}. In these cases, the time step and the damping parameter can be held fixed, confirming the robustness of the proposed approach.

\subsection{Solving 2: Active-set methods (ASM)}
As discussed earlier, standard techniques for contact problems capable to handle box constraints are not straightforward to be applied within the MG iterations, due to the presence of linear constraints. However, a similar approach can be taken in between successive Newton steps. At the Newton step $k$, we check the condition
\begin{equation}
g(x^k)_i \leq \epsilon_x \quad \forall i
\end{equation}
If satisfied, it implies that some of the variable are near to contact the boundary of the feasible region and we deem the corresponding constraints to be active, i.e. $g(x^k)_i=0$, $\forall i\in AS(x^k)$, where $AS(x^k)$ is the set of active constraints. This translates to applying Dirichlet boundary conditions to the variable $\Delta x$ in the systems (\ref{normalEq}) and (\ref{regularizedNormalEq}). This is achieved via the substitutions
\begin{eqnarray}
D_{ii} &=& 1, ~ \forall i \in AS(x^k), \\
A^T_{ij} &=& 0, ~ \forall i \in AS(x^k), \\
(f_x-\nabla g^T G^{-1}f_s)_i &=& -x_k, ~ \forall i \in AS(x^k).
\end{eqnarray}
Moreover, we test for the dual condition
\begin{equation}
\frac{s^k_i}{g(x^k)_i} \leq \epsilon_s.
\end{equation}
This is the deactivation criterion, expressing the fact that if contact forces are not strong enough, then the corresponding constraints are in fact not active. In practice, we found this not to be necessary for the MCF problem and therefore we did not have to deactivate any constraint in the benchmark presented in the next section. With such a treatment of the inequality constraints, our algorithm can be viewed as an active-set method, where rather than testing for a large set of constraints, we use the IPM to find an appropriate guess of the active constraints at the exact solution.

\subsection{Solving 3: Connected components of a graph}
While solutions 1 and 2 are general and can potentially be applied to nonlinear problems, the solution to the third issue exploits the particular nature of the MCF problem. Let us recall the fact that $A^T$ is a discrete gradient operator, it follows that the Lagrange multipliers $y$ are a node-based function whose edge-based gradient is equal to $b$. Therefore, the function $y$ is defined up to a constant. When inequality constraints become active, however, the incidence matrix $A^T$ is modified in a way equivalent to removing the edges corresponding to $AS(x^k)$ from the original graph. Recall the following

\begin{Def}
A weakly connected component is a maximal subgraph of a directed graph such that for every pair of vertices  (u, v) in the subgraph, there is an undirected path from u to v and a directed path from v to u.
\end{Def}

It follows that the multipliers $y$ are defined up to a constant per weakly connected component of the graph. While graph-based algorithms for computing connected components in linear time do exists, we opt instead for taking advantage of parallel matrix-vector multiplication routines. To achieve this, we assemble the adjacency matrix $B$ and perform the following fixed-point iteration:
\begin{eqnarray}
y_i = \max_k B_{ik} y_k,
\end{eqnarray}
where $y$ is initialized as $y_i = i$. The fixed point of this iteration is a vector whose entries sharing the same label belong to the same connected component. Then, we pick one label per connected component and set 
\begin{eqnarray}
L_{ii} &=& 1, ~ \forall j=1..N_{Conn}, \\
L_{iq} &=& 0, ~ \forall j=1..N_{Conn}, ~ q \ne i.
\end{eqnarray}
While not competitive with graph-based algorithms on a single processor, this approach showed to scale very well with the number of cores and is only needed in the close proximity of the true solution, when constraints become active. Hence, further improvements to the run times presented in the next section are possible with a graph-based algorithm, but not to the scaling, which is the primary concern of this work.

\section{Numerical experiments}
We tested the proposed method on the following data sets \cite{Kov15}:
\begin{enumerate}
\item \emph{ROAD\_PATHS}: These instances are road networks from different states in the USA, where the
edge cost is set to the travel time along the edge. There are approximately $\sqrt{n}/10$ randomly selected sources and sinks each with a demand that depends on the maximum flow that can be sent between these sources and sinks.
\item \emph{ROAD\_FLOWS}: Same as above but the capacities are not set to 1, but depending on the category of the
road to either 40, 60, 80 or 100.
\item \emph{NETGEN\_8}: They are generated using the NETGEN random generator. The subscript 8 indicates that the networks are sparse, i.e., $m = 8n$ and the capacities and costs are chosen uniformly at random between 1 and 1000 and 1 and 10000, respectively. There are approximately $\sqrt{n}$ sinks and sources and the total demand is $1000 \sqrt{n}$.
\item \emph{GOTO\_8}: They are generated using the GOTO generator. They are grid instances on tori, known to be rather hard. The subscript 8 indicates that the networks are sparse, i.e., $m = 8n$ and the capacities and costs are chosen uniformly at random between 1 and 1000 and 1 and 10000, respectively. There is one source and one sink node and the supply is adjusted to the capacity.
\item \emph{Signal processing}: These instances are regular quadrilateral grids taken from satellite data that can be measured only up to mod $2\pi$. Costs and flows are specified in order to reconstruct the true values of the function over the grid \cite{Cos98}.
\end{enumerate}
For our implementation, we rely on BICGStab from PETSc \cite{Bal15} as the linear solver and the BoomerAMG algebraic MG preconditioner, which is part of the Hypre package \cite{Fal02}. The default Hypre options and parameters have been used in all data sets, with the exception of NETGEN and GOTO, for which we discuss below how AMG has been tuned. 

\subsection{Summary of results} We are interested in three types of scaling:

\begin{itemize}
\item Scaling with respect to problem size. The results are summarized in Figure \ref{exponent_table}, using the exponent $\alpha$ of the best fit of run times to the complexity model $m^\alpha$, where $m$ is the total number of arcs. As it can be seen, the scaling of the proposed algorithm is better than the best known MCF algorithms in all cases except for the NETGEN data set. However, as discussed below, there is a moderately large constant complexity factor ($\approx 200$) with our method compared to combinatorial approaches that is not captured by this estimate, making them more efficient for a small problem with a small number of cores. Ad discussed below, such a constant is heavily dependent on the choice of the MG method and might be improved by using a graph Laplacian AMG.

\begin{figure}[ht!]  \centering 
\fbox{
\begin{tabular}{c|c|c|c|c}
Data set & Ours & SSP & NS & COS \\ \hline
R\_PATHS &  \textbf{1.3371}  & 1.4065 & 1.7234 & 1.4274   \\
R\_FLOWS &  \textbf{1.3399} & 1.4292 & 1.7212 & 1.4423   \\
NETGEN & 2.0797 & 1.8860 & 1.6824 & \textbf{1.2205} \\
GOTO & \textbf{1.4725} &  2.1076 & 2.2676 & 1.5356   \\
Signal proc. & \textbf{1.0875} & - & - & -
\end{tabular}
}
\caption{Summary of scaling exponents for our method against the network simplex algorithm (NS), the cost scaling algorithm (COS), and the successive shortest path method (SSP). Scaling for our method was computed on two cores on a laptop for the first four data sets and on two nodes (48 cores) for the signal processing case.}\label{exponent_table}
\end{figure} 

\item Scaling with respect to problem complexity. In this case, an exact measure is difficult to be defined, but it is useful to compare the results on the different benchmarks, as done in Figure \ref{PETSc_table}, where a breakdown of run times for the most expensive routines is given for our algorithm. From this comparison, it is clear that NETGEN and GOTO are particularly difficult to be solved. However, such a poor performance can be mostly explained by BoomerAMG, which had to be tuned to avoid producing dense interpolation operators exceeding the memory availability, at the expense of performance. On the other hand, on more regular grids, combinatorial algorithms proved to be much more sensitive to problem complexity, taking 650\% more time when solving ROAD\_FLOW compared to ROAD\_PATHS, while our method showed an increase of only 50\%. 

\begin{figure}[ht!]  \centering 
\fbox{
\begin{tabular}{c|c|c|c|c|c|c|c}
Data set & Arcs & Cores & Total & BCGS & PCSet & PCAppl & MatMult  \\ \hline
R\_PATHS &  5.2m    & 192 & 29s & 76\% & 10\% & 62\% & 14\% \\
R\_FLOWS &  5.2m &192 & 47s & 77\% & 11\% & 64\% & 13\% \\
NETGEN & 33.5m &192 & 4721s & 98\% & 39\% & 56\% & 3\% \\
GOTO & 4.2m &96  & 3556s & 43\% & 27\% & 16\% & 47\% \\
Signal proc. & 100m & 300  & 580s & 42\% & 6\% & 33\% & 34\%
\end{tabular}
}
\caption{Breakdown of computational effort for the largest instance in each data set. PCSetup and PCApply refer to the routines from BoomerAMG. }\label{PETSc_table}
\end{figure} 

\item Weak scaling tests with respect to number of cores and fixed data set, where the number of arcs per core is kept constant as the problem size is increased. Given that the expected complexity of IPM is $O(n^{1.5})$, we define the weak scaling efficiency $\eta_k$ of the $k$-th instance as
\begin{eqnarray} \label{Efficiency}
\eta_k = 100 \left( \frac{m_k}{m_0} \right)^{1.5}\frac{\#\mathrm{cores}_0}{\#\mathrm{cores}_k}\frac{\mathrm{time}_0}{\mathrm{time}_k}.
\end{eqnarray}
The results concerning weak scaling are given separately for each data set in the remainder of this section.
\end{itemize}

\subsection{Graph Laplacian AMG}
As already discussed, from Figure \ref{exponent_table} it is clear that the performance on the NETGEN data set is very poor. In order to confirm that this is a limitation of BoomerAMG, we have tested our approach with the LAMG method \cite{Liv12}, which is an algebraic MG solver specifically made for the graph Laplacian, for which a Matlab implementation is available. We have obtained a nearly-optimal scaling factor of $1.0318$, better than any known combinatorial algorithm, comparable to the BoomerAMG performance on a regular grid network, such as the signal processing data set. This gives a good evidence that an optimal-complexity IPM solver for the MCF on massively parallel architectures is within reach with our strategy, once an efficient parallel implementation of LAMG becomes available.  At the same time, we have also tested the application of LAMG to the IPM without our proposed modification and found that both BoomerAMG and LAMG fail to converge near the solution. This agrees with our discussion in Section 3 and confirms the need to regularize the problem in order to apply efficient solvers.

\subsection{Weak scaling tests}
We present here the results from weak scaling tests. The hardware used is the Piz Dora supercomputer at the CSCS center. It is a Cray XC40 system with 1256 compute nodes, each equipped with 64 GB of RAM and two 12-core Intel Haswell CPUs (Intel Xeon E5-2690 v3). 

\subsubsection{ROAD\_PATHS}
\begin{figure}[!ht]
\begin{center}
\includegraphics[width=.9\textwidth]{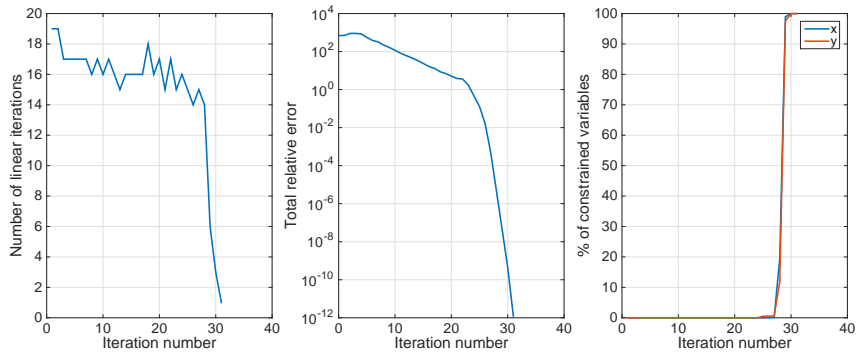}
\caption{Number of linear solver iterations, total relative error, and percentage of active variables (primal and dual), as a function of the IPM iteration number for the road paths instance TX\_07. }
\label{road_paths_plots}
\end{center}
\end{figure}

As presented in \cite{Kov15}, the fastest combinatorial algorithm takes 21.99 seconds on the largest instance (TX\_07) of this data set. The results for our method are shown in Figures \ref{road_paths_table} and \ref{road_paths_plots}. As it can be seen, our strategy requires about 200 cores to match the single-core performance of the best combinatorial algorithm. However, despite the fact that 25872 unknowns per core can be considered too small as a load distribution, the efficiency is above optimal. As can be seen in Figure \ref{exponent_table}, this is also due to the fact that the method scales with $O(m^{1.3371})$ rather than $O(m^{1.5})$. Therefore, the need for 200 cores is not a symptom of poor scaling but it is a constant complexity factor, so our algorithm would be able to efficiently tackle a larger instance of this data set. From Figure \ref{road_paths_plots}, it becomes clear that our modifications to the problem improve the convergence behavior of the linear solver. It is worth remarking that had they been disabled, then rather than decrease, the number of linear iterations would blow up in the last Newton iterations, until the point where BoomerAMG would completely fail to converge.

\begin{figure}[ht!]  \centering 
\fbox{
\begin{tabular}{c|c|c|c|c|c}
Data set & Cores (\#) & Time (s) & Efficiency (\%) & Arcs per core & It (\#)  \\ \hline
04 NV &      24 & 12.7 & 100 & 25872 & 29 \\
05 WI & 48 & 13.2 & 136 & 25872 & 25 \\
06 FL & 96 & 19.4 & 131 & 25872 & 27 \\
07 TX & 192  & 29.2 & 123 & 25872 & 31
\end{tabular}
}
\caption{Summary of road paths results. The last column refers to the total number of Newton iterations. }\label{road_paths_table}
\end{figure} 

\subsubsection{ROAD\_FLOW}
\begin{figure}[!ht]
\begin{center}
\includegraphics[width=.9\textwidth]{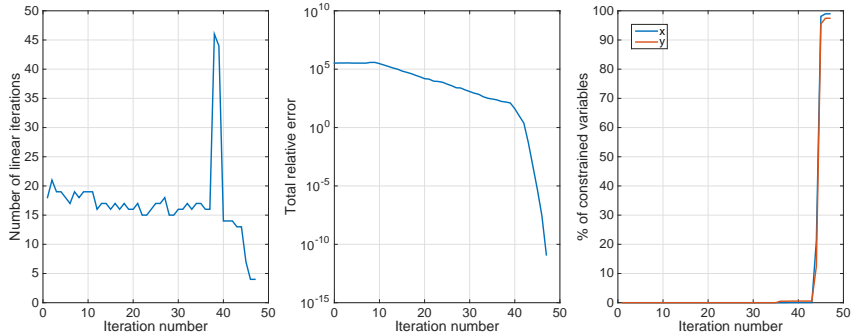}
\caption{Number of linear solver iterations, total relative error, and percentage of active variables (primal and dual), as a function of the IPM iteration number for the road paths instance TX\_07. }
\label{road_flows_plots}
\end{center}
\end{figure}

As presented in \cite{Kov15}, the fastest combinatorial algorithm takes 145.7 seconds on the largest instance (TX\_07). The results for our  method are shown in Figures \ref{road_flows_table} and \ref{road_flows_plots}. By comparing this table with the ROAD\_PATHS scenario, we see that as the complexity of the problem increases, run times increase by only $50\%$ with our method, while they increase by more than $650\%$ with combinatorial algorithms. It can also be seen that efficiency is lower in this case, due to a larger increase in the number of Newton iterations required to achieve convergence. We can also notice from Figure \ref{road_flows_plots}, that the number of linear iterations has a spike around the Newton iteration 38. We argue that this is an indication that an adaptive selection of the $\rho$ and $\Delta t$ parameters could give a performance benefit.

\begin{figure}[ht!]  \centering 
\fbox{
\begin{tabular}{c|c|c|c|c|c}
Data set & Cores (\#) & Time (s) & Efficiency (\%) & Arcs per core & It (\#) \\ \hline
04 NV &      24 & 13.2 & 100 & 25872 & 31 \\
05 WI & 48 & 17.9 & 104 & 25872 & 35  \\
06 FL & 96 & 42.5 & 62 & 25872 & 40 \\
07 TX & 192  & 47 & 79 & 25872 & 48 
\end{tabular}
}
\caption{Summary of road flows results. The last column refers to the total number of Newton iterations.}\label{road_flows_table}
\end{figure} 

\subsubsection{NETGEN}
\begin{figure}[!ht]
\begin{center}
\includegraphics[width=.9\textwidth]{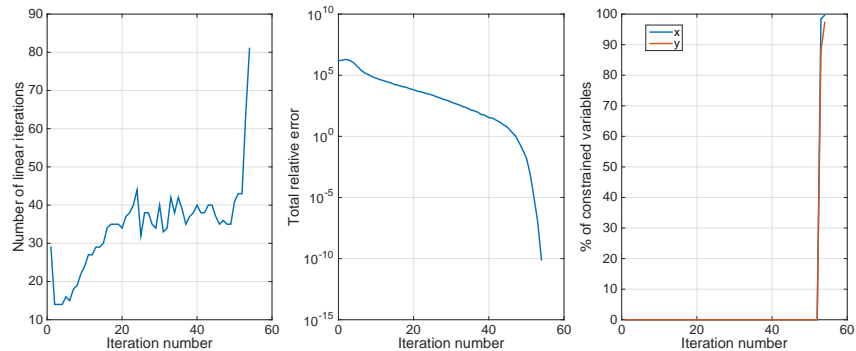}
\caption{Number of linear solver iterations, total relative error, and percentage of active variables (primal and dual), as a function of the IPM iteration number for the NETGEN instance 22a. }
\label{netgen_plots}
\end{center}
\end{figure}

As presented in \cite{Kov15}, the fastest combinatorial algorithm takes 419.4 seconds on the largest instance (22a). The results for the proposed method are shown in Figures \ref{netgen_table} and \ref{netgen_plots}. For this data set, we have observed that Newton and MG both converge in a reasonable amount of iterations, however BoomerAMG is struggling in the setup phase, resulting is a very slow assembly of the coarse problems and in a very high memory occupation. We have limited these issues by using 5 levels of aggressive coarsening, at the price of an increased number of linear iterations, as can be seen from Figures \ref{PETSc_table} and \ref{netgen_plots}. Unfortunately, this is an intrinsic limitation of algebraic MG methods and the implementation of a custom geometric MG would be required to considerably improve efficiency. This will be the subject of future studies.

\begin{figure}[ht!]  \centering 
\fbox{
\begin{tabular}{c|c|c|c|c|c}
Data set & Cores (\#) & Time (s) & Efficiency (\%) & Arcs per core & It (\#) \\ \hline
19a & 24 & 900 & 100 & 174763& 42 \\
20a & 48 & 1445 & 88 & 174763 & 46 \\
21a & 96  & 2327 & 77 & 174763 & 50 \\
22a & 192  & 4721 & 54  & 174763 & 54
\end{tabular}
}
\caption{Summary of NETGEN results. The last column refers to the total number of Newton iterations.}\label{netgen_table}
\end{figure} 

\subsubsection{GOTO}
\begin{figure}[!ht]
\begin{center}
\includegraphics[width=.9\textwidth]{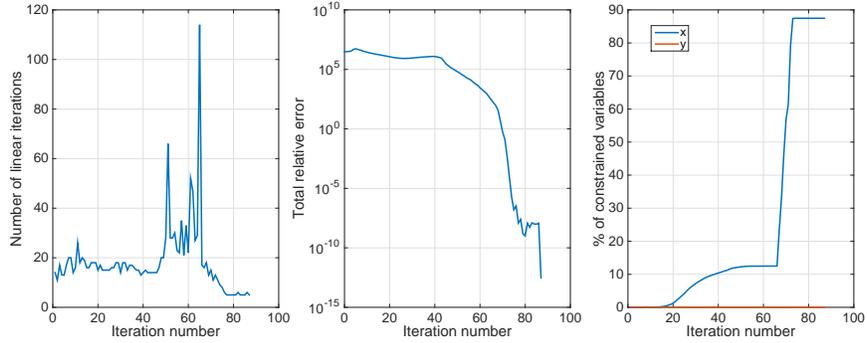}
\caption{Number of linear solver iterations, total relative error, and percentage of active variables (primal and dual), as a function of the IPM iteration number for the GOTO instance 19a. }
\label{goto_plots}
\end{center}
\end{figure}

As presented in \cite{Kov15}, the fastest combinatorial algorithm takes 1564 seconds on the second largest instance (19a), while no data is given for the largest one (22a). It is worth remarking that this is a hard problem and run times for most of the combinatorial algorithms are not given beyond the 16a instance. The results for the proposed method are shown in Figures \ref{goto_table} and \ref{goto_plots}. For this example, we have observed that the number of Newton iterations grow considerably as the instance size increases. To obtain a reasonable performance, we used the adaptive time step strategy presented in Algorithm \ref{adaptiveDt} with $\rho=\epsilon=10^{-6}$. Despite this, it is clear from the first plot in Figure \ref{goto_plots} that the number of linear iterations is far from optimal. This is also due to the fact that, as done in the NETGEN case, we had to limit the memory and computational effort of the BoomerAMG setup phase, by restricting the number of nonzero in the interpolation operator to 5. Therefore, it is possible that also in this case a different MG implementation would improve the convergence behavior. As can be seen from the last plot in Figure \ref{goto_plots}, this is the only example where the active set becomes already active in the first iterations. In general, this behavior is not ideal, but we found it to be necessary here in order to achieve MG convergence.

\begin{figure}[ht!]  \centering 
\fbox{
\begin{tabular}{c|c|c|c|c|c}
Data set & Cores (\#) & Time (s) & Efficiency (\%) & Arcs per core & It (\#)  \\ \hline
16a & 12 & 118 &  100 & 43690 & 56 \\
17a & 24 & 332 & 50  & 43690 & 72 \\
18a & 48  &   820 & 29  & 43690 & 51 \\
19a & 96 & 3556 &  9 & 43690 & 95 
\end{tabular}
}
\caption{Summary of GOTO results. The last column refers to the total number of Newton iterations.}\label{goto_table}
\end{figure} 

\subsubsection{Signal processing}
The results for the proposed method are shown in Figure \ref{signal_table}. If we compare Figures \ref{signal_table} and \ref{netgen_table}, we immediately see that BoomerAMG is 8x faster on this example. This is not surprising, since the regularity of the mesh makes this example much more similar to FE applications, for which BoomerAMG was specifically developed. This also suggests that there would be a considerable speedup available in the NETGEN case from implementing a custom MG solver for graph Laplacians. However, this example shows that our method is well-suited to treat problem sizes that are difficult to solve in a reasonable time by a combinatorial algorithm on a single core.

\begin{figure}[ht!]  \centering 
\fbox{
\begin{tabular}{c|c|c|c|c|c}
Data set & Cores (\#) & Time (s) & Efficiency (\%) & Arcs per core & It (\#)  \\ \hline
1000x1000 & 12 & 150 &  100 & 333333 &  33 \\
2000x2000 & 48 & 206 &  146  & 333333 & 33 \\
3000x3000 & 108 & 219 &  205 & 333333 &  31 \\
4000x4000 & 192 & 312 &  192 & 333333 &  32 \\
5000x5000 & 300 & 580 &  129 & 333333 &  45
\end{tabular}
}
\caption{Summary of signal processing results. The last column refers to the total number of Newton iterations.}\label{signal_table}
\end{figure} 

\begin{figure}[!ht]
\begin{center}
\includegraphics[width=.9\textwidth]{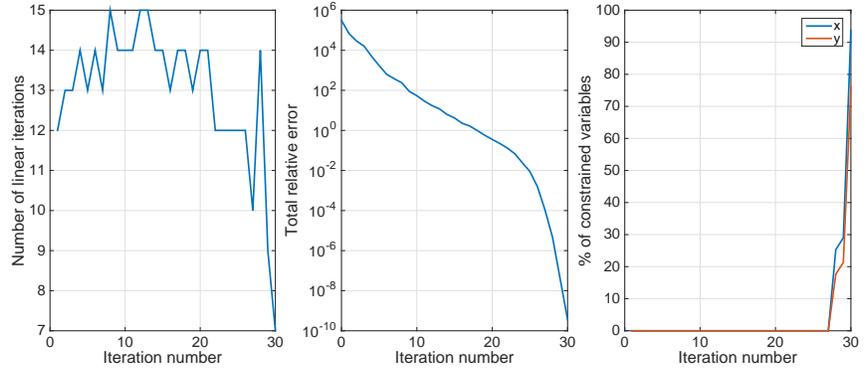}
\caption{Number of linear solver iterations, total relative error, and percentage of active variables (primal and dual), as a function of the IPM iteration number for the signal processing instance 3000x3000. }
\label{meshes}
\end{center}
\end{figure}

\section{Conclusions}
We have presented an algorithm to solve the MCF problem on massively parallel architectures. Our approach is based on interpreting the Newton iteration arising from solving the KKT equations as a mixed formulation of the graph Laplacian of the network, with an appropriate metric term that depends on the inequality constraints slack and multipliers. We have proposed an intertwining of IPM and ASM in order to prevent such a metric from vanishing or blowing up, so that it is possible to apply optimal-complexity linear solvers for elliptic operators, such as AMG. We have shown with standard benchmarks that, while less competitive than combinatorial techniques on small problems that run on a single core, our approach scales well with problem size, complexity, and number of processors, allowing for tackling large problems in a reasonable time.

There are future improvements that can be used to enhance our strategy:
\begin{enumerate}
\item A geometric MG or a parallel graph Laplacian AMG can be used to replace the BoomerAMG solver. This would allow to take advantage of the geometric or graph informations in order to select a more appropriate grid hierarchy for the metric matrix $D$, greatly improving efficiency and memory consumption. Indeed, we have shown that nearly optimal-complexity is already achievable when the network geometry is regular, such as for the signal processing data, or when the LAMG solver made for graph Laplacians is used. This gives a good indication that an optimal-complexity parallel solver for the MCF problem is within reach, once a parallel implementation of LAMG is available.
\item An optimal choice of the parameters $\mu$, $\rho$, $\beta$, and $\Delta t$, can be studied for achieving optimal damping ratio of 1, in order to reach the static equilibrium as quickly as possible, and to ensure that the total sum of linear solver iterations over the Newton steps is minimized. 
\end{enumerate}


\bibliographystyle{siamplain}

\end{document}